\documentclass[12pt]{article}
\usepackage{amsmath,amssymb}
\newcommand{\G}{\Gamma}
\newcommand{\nin}{\noindent}
\newcommand{\vs}{\vspace*}
\newcommand{\al}{\alpha}

\newcommand{\seq}{\subseteq}

\newcommand{\ol}{\overline}

\title {\bf The structures of a class of $Z$-local rings }
\author{Tongsuo Wu\thanks{Corresponding author. {\small email:
tswu@sjtu.edu.cn}}\\
Department of Mathematics\\  Shanghai Jiaotong University\\
Shanghai 200030, P. R. China\\ and \\Dancheng Lu\thanks{\small
email: ludancheng@suda.edu.cn}\\
Department of Mathematics\\  Suzhou University\\
Suzhou 215006, P. R. China\\ }
\date{}
\begin{document}
\baselineskip=16pt \maketitle \vspace{4mm}
\begin{center}
\begin{minipage}{13cm}
\vs{3mm}\nin{\bf Abstract:} {\small A local ring $R$ is called
$Z$-local if $J(R)=Z(R)$ and $J(R)^2=0$. In this paper the
structures of a class of $Z$-local rings are determined.}

\vs{2mm}\nin {\small\bf Key Words:} {\small $Z$-local ring,
structure, polynomial rings}
\end{minipage}
\end{center}

\vs{3mm} Let $R$ be a commutative local ring which is not
necessarily noetherian. Denote by $J(R)$ the Jacobson radical of
$R$, $Z(R)$ the zero-divisor elements of $R$. $R$ is called {\it
$Z$-local} if $J(R)=Z(R)$ and $J(R)^2=0$. This concept was
introduced in [2] where the authors proved that for any
commutative ring $S$ such that $2$ is regular in $S$ and that $S$
satisfies DCC on principle ideals, if the zero-divisor graph
$\G(S)$ of $S$ is uniquely determined by neighborhoods and $S$ is
not a Boolean ring, then $S$ is a $Z$-local ring. The zero-divisor
graph of a commutative ring was introduced and studied in [1]. In
this paper, we will try to determine the structure of a class of
$Z$-local rings.

For any commutative ring extension $A\subseteq  B$ and any $\al\in
B$, recall that $\al$ is said to be {\it integral over $A$}, if
there is a monic polynomial $f(x)\in A[x]$ such that $f(\al)=0$. It
is well known that $\al$ is integral over $A$ if and only if there
is a subring $C$ of $B$ which contains $A$, such that $\al\in C$ and
{\it $C$ is finitely generated as an $A$-module} (Please see, e.g.,
[3, Theorem 9.1]). For an element $\al$ integral over $A$, a minimal
polynomial of $\al$ over $A$ is a monic polynomial $f(x)$ with the
least degree such that $f(\al)=0$. In general, a minimal polynomial
over $A$ need not be unique. But if $A$ is a field, then it is
unique. If $A=\Bbb Z_{p^2}$ for some prime number $p$, then {\it the
minimal polynomial of $\al$ over $\Bbb Z_{p^2}$ is unique modulo
$p$}, that is, $p$ divides all the coefficients of $f(x)-g(x)$ for
any minimal polynomials $f(x)$ and $g(x)$ of $\al$ over $\Bbb
Z_{p^2}$. In this case, we will denote it as $$f(x)\equiv g(x)(\mod
p).$$ \nin These observations will be used in the latter part of the
paper.

By [2,Theorem 2.5], the characteristic of a $Z$-local ring has
only three possible values, i.e., 0, $p$ or $p^2$. For a $Z$-local
ring $R$ with characteristic $0$, since any element of a $Z$-local
ring is either a unit or a zero-divisor, we have $\Bbb Q\seq
T(R)\cong R$.

\vs{3mm}\nin{\bf Theorem 1.} {\it For a $Z$-local ring $R$, let
$F$ be the prime subfield of the field $K=R/J(R)$. Assume that the
characteristic of the ring $R$ is not $p^2$ for any prime number
$p$. Assume further that $K=F[\ol{\alpha}]$ is an algebraic
extension over $F$ for some $\alpha\in R$, and let $g(x)$ be the
minimal polynomial of $\ol{\alpha}$ over $F$ with degree $n$. Let
$<S>=\{s_i\,|\,i\in I\}$ be the $K$-basis of the $K$-module
$J(R)$.  }

(1){\it There is an $F$-algebra epimorphism from $F[ x,Y]/<\{
y_iy_j\,|\,i,j\in I\}>$ to $R$, where $Y=\{y_i\,|\,i\in I\}$ is a
set of commutative indeterminants.}

(2) {\it If $\al$ is integral over F, then one and only one
situation occurs in the following:}

($i$) {If $g(\alpha)=0$, then $R\cong K[Y]/<\{ y_iy_j\,|\,i,j\in
I\}>,$ where $g(x)$ is irreducible over $F$ and $m\ge 1$.}

($ii$) {\it If $g(\al)\not=0$, then assume
$g(\al)=\sum_{i=1}^mv_i(\al)s_i$, where $v_i(\al)\not\in J(R)$.
Then
$$R\cong F[x,Y]/<g(x)^2,g(x)-\sum_{i=1}^mv_i(x)y_i,g(x)y_r,y_sy_t|r,s,t\in I>.$$
(Notice that in the second case, $m\ge 2$, $g(x)$ is irreducible
over $F$ with degree at least 2, and $v_i(x)$ are nonzero
polynomials over $F$ and $deg(v_i(x))<n$. )}

\vs{3mm}\nin{\bf Proof.} {\bf (1).} By assumption, we have
$F\subseteq R$. Since
$$(F[\alpha]+J(R))/J(R)=F[\ol{\alpha}]=R/J(R),$$
\nin we have $R=F[\alpha]+J(R)$. Now consider the $F$-algebra
homomorphism
$$\sigma:F[x,Y]/<\{
y_iy_j\,|\,i,j\in I\}>\to R=F[\alpha]+J(R),\ol{h(x,y_i)}\mapsto
h(\alpha, s_i ).$$ \nin By assumption $S$ is the set of generators
of the $R$-module $J(R)$. Since $J(R)^2=0, R=F[\alpha]+J(R)$, thus
$\sigma$ is a surjective $F$-algebra homomorphism. This proves the
first part of the theorem.

{\bf (2).} Now assume further that $\al$ is integral over F. Let
$g(x)\in F[x]$ be the minimal monic polynomial of $\ol{\alpha}$
over $F$ and assume $deg(g(x))=n$. Let $f(x)$ be the minimal monic
polynomial of $\alpha$ over $F$. Then $g(x)$ is irreducible in
$F[x]$ and we have $g(x)|f(x)$. Now assume $f(x)=g(x)^u\cdot
l(x)$, where $(l(x),g(x))=1$. Since $g(\al)\in J(R)$, thus
$l(\al)\in U(R)$. By assumption. we must have $f(x)=g(x)^u$, where
$u\le 2$.

{\bf Case 1.} If $u=1$, then in $R$ we have $g(\al)=0$. By
assumption, for each nonzero polynomial $r(x)$ of degree less than
$n$ ($n=\mbox{deg}((g(x)) $), we have $r(\al)\not\in J(R)$, i.e.,
$r(\al)$ is invertible in $R$. Thus $R=F[\al]\oplus J(R)$, where
$F[\al] $ is a field  and it is also a subring of $R$. In this
case, we obviously have an $F$-algebra isomorphism
$$R\cong K[Y]/<\{ y_iy_j\,|\,i,j\in I\}>.$$

{\bf Case 2.} If $u= 2$, then $g(\al)\not=0$ and $g(\al)\in J(R)$.
\nin In this case, consider
$$\tau:F[x,Y]/W\to R, \ol{h(x,y_i)}\mapsto h(\al,s_i ),$$
\nin where
$$W=<g(x)^2,g(x)-\sum_{i=1}^mv_i(x)y_i,g(x)y_r,y_sy_t|r,s,t\in I>.$$
\nin By assumption, $\tau$ is a map and thus a surjective
$F$-algebra homomorphism.In order to prove that $\tau$ is
injective, for any $h(x,y_i)\in F[x,Y]$, we have
$$h(x,y_i)=g(x)^2A+\sum_{i,j}y_iy_jB_{ij}+(\sum_{r}[g(x)q_r(x)+f_r(x)]y_r+[g(x)q_{\#}(x)+f_{\#}(x)],$$
\nin where the degrees of $q_{d}(x)$ and $f_{e}(x)$ are at most
$n-1$ whenever they are not zero. By assumption, $q_d(\al)$ and
$f_{e}(\al)$ are units of $R$ when the corresponding polynomials
are not zero. Then if $h(\al, s_i)=0$, then we must have
$$0=(\sum_{r}[f_r(\al)]s_r+[g(\al)q_{\#}(\al)+f_{\#}(\al)],\quad (*)$$
Now if $f_{\#}(x)\not=0$, then $f_{\#}(\al)$ is a unit. But from
the previous equality ($*$), we also obtain $f_{\#}(\al)\in J(R)$,
a contradiction. Thus by assumption and $(*)$, we obtain
$$\sum_{r}f_r(\al)\cdot s_r=- q_{\#}(\al)\sum_{i=1}^m v_i(\al)s_i.$$
Thus for $r\not\in \{1,2,\cdots,m\}$, $f_r(x)$ must be zero or
else it is a unit and at the same time, it is in $J(R)$. For
$r=1,2,\cdots,m$, we have $f_r(\al)=-q_{\#}(\al)\cdot v_r(\al)$.
Since $f(x)=g(x)^2$, we obtain $f_r(x)=-q_{\#}(x)\cdot v_r(x).$
Now coming back to the previous decomposition of $h(x,y_i)$, we
obtain
$$\sum_{r}f_r(x)\cdot y_r+g(x)q_{\#}(x)=q_{\#}(x)[g(x)-\sum_{i=1}^mv_i(x)y_i].$$
\nin This shows that $\tau$ is injective. This completes the whole
proof. $\Box$

\vs{3mm}Now let $R$ be a $Z$-local ring with $\mbox{char}(R)=p^2$
for some prime number $p$. Then $\{i\,|\,0\le i\le  p^2-1\}=\Bbb
Z_{p^2}\subseteq R$. Denote by $F$ the prime subfield $\Bbb Z_p$ of
the field $K=R/J(R)$ and let $S\cup \{p\}$ be a set of $K$-basis of
the $K$-space $J(R)$, where $p\not\in S$ and $S=\{s_i\,|\,i\in I\}$.
Let $Y=\{y_i\,|\,i\in I\}$ be a set of indeterminants determined by
the index set $I$. Assume further that $K=F[\ol{\alpha}]$ is an
algebraic extension over $F$ for some $\alpha\in R$, and let
$\ol{g}(x)$ be the minimal polynomial of $\ol{\alpha}$ over $F$ with
degree $n$, where $g(x)\in \Bbb Z_{p^2}[x]$ is a monic polynomial.
We also observe the following facts which will be used repeatedly:

{\bf For any polynomial $u(x)\in \Bbb Z_{p^2}[x]$, if
$u(x)\not\equiv 0 (\mod p)$ and its degree modulo $p$ is less than
n, then $u(\al)$ is a unit of $R$.}

We are now ready to determine the structure of a class of
$Z$-local rings with characteristic $p^2$.

 \vs{3mm}\nin{\bf Theorem 2.} {\it For a $Z$-local ring $R$ with characteristic $p^2$,
assume that $R/J(R)=K=F[\ol{\alpha}]$ is an algebraic extension over
$\Bbb Z_p\cong F\seq R/J(R)$ for some $\alpha\in R$. Assume further
that $\al$ is integral over $\Bbb Z_{p^2}$. Then either $R\cong \Bbb
Z_{p^2}[x,Y]/Q_1$, where $Q_1=<Q\cup\{g(x)\}>$ and $|Y|\ge 0$, or
$R\cong \Bbb Z_{p^2}[x,Y]/Q_2$, where
$$Q_2=<Q\cup\{g(x)^2, pg(x),g(x)y_r,
g(x)-\sum_{i=1}^{m}v_i(x)y_{i}\,|\,r\in I\}>.$$ \nin  In each case,
$\ol{g}(x)$ is irreducible over $\Bbb Z_p$, and
$$Q=\{px, y_sy_t,py_r,|r,s,t\in
I >,$$ \nin where $m\ge 1$ is a fixed number, and $v_i(x)\not\equiv
0(\mod p)$. We also notice that in the second case, $g(x)$ is some
polynomial over $\Bbb Z_{p^2}$ such that} $\mbox{deg}$\,$g(x)>1$.

\vs{3mm}\nin{\bf Proof.} First, it is easy to see that $R=\Bbb
Z_{p^2}[\al]+J(R)$. Since $\al$ is integral over $\Bbb Z_{p^2}$, we
have a minimal polynomial $f(x)\in \Bbb Z_{p^2}[x]$ which is unique
modulo $p$. By the choice of $g(x)$, we have $f(x)=g(x)^ur(x)(\mod
p)$ for some monic $r(x)\in \Bbb Z_{p^2}[x]$ satisfying $(\ol{g}(x),
\ol{r}(x))=\ol{1}$ in $F[x]$. Thus $r(\al)$ is invertible in $R$
since $g(\al)\in J(R)$. Without loss of generality, we can assume
$f(x)=g(x)^ur(x)+h(x)$, where $h(x)\equiv 0(\mod p)$. If $u\ge 3$,
then we obtain $0=f(\al)=g(\al)^ur(\al)+h(\al)=h(\al)$. Thus the
monic polynomial $g(x)^2r(x)$ annilates $\al$ and it has a degree
less than deg$f(x)$, contradicting to the choice of $f(x)$. Thus we
must have $u\le 2$.

{\bf Case 1.} If $u=2$, we have $h(\al)=0$ again. In this case, we
must have $r(x)=1$ since $g(x)^2$ annilates $\al$. In this case,
we can choose $f(x)=g(x)^2$.

{\bf Case 2.} If $u=1$, we have $f(x)=g(x)r(x)+h(x)$, where
$h(x)\equiv 0(\mod p)$. In this case, we have $g(\al)=-h(\al)\cdot
r(\al)^{-1}=-h(\al)w(\al)$ for some $w(x)\in \Bbb Z_{p^2}[x]$.
Obviously $g(x)\equiv g(x)+h(x)w(x)(\mod p)$. Thus in this case we
can choose the $g(x)$ such that $g(\al)=0$.

\vs{3mm}{\bf (1)} Let us first consider the case when $g(\al)=0$.

 In this case, consider
$$\tau:F[x,Y]/Q_1\to R=F[\alpha]+J(R),\ol{h(x,y_i)}\mapsto
h(\alpha, s_i ).$$ \nin For each generators $h(x,y_i)$ of $Q_1$,
we have $h(\alpha, s_i )=0$. Thus $\tau$ is a surjective
$F$-algebra homomorphism. Now for any $h(x,y_i)\in F[x,Y]$, we
have a decomposition
$$h(x,y_i)\equiv \sum_{r}f_r(x)y_r+f_{\#}(x)(\mod Q_1),\quad (**)$$
\nin where $f_r(x)$ are some polynomials of $x$ over $F$ which has
degree less than $n$ when they are nonzero modulo $p$, for all
$s\not=\#$. If in $F[x,Y]/Q_1$, $\ol{h(x,y_i)}\not=0$, then either
one of the $f_s(x)$ is not zero modulo $p$, or $f_{\#}(x)\not=0$.
Thus if $f_{\#}(x)=0$, then we have some unit $f_s(\al)$ and hence
$h(\al,s_i)\not=0$. If $f_{\#}(x)\not=0$, we also conclude that
$h(\al,s_i)\not=0$. In fact, assume in the contrary that
$h(\al,s_i)=0$. If deg$(f_{\#}(x))>0$ with coefficients modulo
$p$, then $f_{\#}(\al)\in J(R)\cap U(R)$, a contradiction. If
$f_{\#}(x)\not=0$ and deg$(f_{\#}(x))=0$, then we need only
consider the case when $f_{\#}(\al)=pi (\mod px)$ for some $1\le
i\le p-1$ since $px\in Q_1$. Then we obtain a contradiction
$i\cdot p+\sum_{r}f_r(\alpha)s_r=0$, since $p\not\in S$. These
arguments show that $\tau$ is injective. In conclusion, $\tau$ is
an $F$-algebra isomorphism under the assumption of $g(\al)=0$.

{\bf (2)} Now assume $g(\al)\not=0$. Then $g(\al)^2=0$ and
$pg(\al)=0$ since $g(\al)\in J(R)$. Since this case corresponds to
the case of $f(x)=g(x)^2$, we can choose an $g(x)$ such that
$g(\al)=\sum_{i=1}^m v_i(\al)s_i$, where $v_i(\al)\in U(R)$.

In this case, consider
$$\tau:F[x,Y]/Q_2\to R, \ol{h(x,y_i)}\mapsto h(\al,s_i ),$$
\nin By the choice of $Q_2$, it is easy to see that $\tau$ is a
map and thus a surjective $F$-algebra homomorphism. In order to
prove that $\tau$ is injective, for any $h(x,y_i)\in F[x,Y]$, we
have
$$h(x,y_i)\equiv \sum_{r}f_r(x)y_r+[g(x)q_{\#}(x)+f_{\#}(x)](\mod Q_2),$$
\nin where the degrees of $q_{\#}(x)$ and $f_{r}(x)$ are at most
$n-1$ whenever they are not zero, modulo $p$. By assumption,
$q_{\#}(\al)$ and $f_{r}(\al)$ are units of $R$ when the
corresponding polynomials are not zero modulo $p$ ($r\not=\#$). If
$h(\al, s_i)=0$, then we must have
$$0=(\sum_{r}[f_r(\al)]s_r+[g(\al)q_{\#}(\al)+f_{\#}(\al)],\quad (*)$$

{\bf (Subcase 1.)} If $f_{\#}(x)=0$, then by assumption and $(*)$,
we obtain
$$\sum_{r}f_r(\al)\cdot s_r=- q_{\#}(\al)\sum_{i=1}^m v_i(\al)s_i.$$
Thus for $r\not\in \{1,2,\cdots,m\}$, $f_r(x)$ must be zero
(modulo $p$), or else $f_r(\al)$ is a unit and at the same time,
it is in $J(R)$. For $r=1,2,\cdots,m$, we have
$f_r(\al)=-q_{\#}(\al)\cdot v_r(\al)$ ($\mod J(R)$). Since
$f(x)=g(x)^2$, we obtain $f_r(x)=-q_{\#}(x)\cdot v_r(x)$ modulo
$p$. Now coming back to the previous decomposition of $h(x,y_i)$,
we obtain
$$\sum_{r}f_r(x)\cdot y_r+g(x)q_{\#}(x)=q_{\#}(x)[g(x)-\sum_{i=1}^mv_i(x)y_i]\equiv 0 (\mod Q_2).$$

{\bf (Subcase 2.)} If $f_{\#}(x)\not=0$, then deg$(f_{\#}(x))=0$
(modulo p), or else $f_{\#}(\al)\in J(R)\cap U(R)$ by $(*)$, a
contradiction. In the following, we assume $f_{\#}(x)\not=0$ and
deg$(f_{\#}(x))=0$ (with coefficients  modulo $p$). Now consider
$(*)$ and assume $f_{\#}(x)\equiv pi(\mod px)$ for some $1\le i\le
p-1$. We have
$$i\cdot p+ \sum_{r}f_r(\al)\cdot s_r\equiv -q_{\#}(\al)g(\al).$$
\nin Since $i$ is invertible in $R$, $p$ can be written as an
$R$-combination of the $s_i$'s. This is certainly impossible. The
above arguments show that $\tau$ is injective. This completes the
whole proof. $\Box$

\vs{3mm}It is well known that any finite field is a simple algebraic
extension over it's prime subfield $\Bbb Z_p$ for some prime number
$p$. As an application of Theorem 2, we immediately obtain the
following results.

\vs{3mm}\nin{\bf Theorem 3.} {\it Let $R$ be a finite ring whose
characteristic is a prime square $p^2$. Then $R$ is a $Z$-local
ring if and only if either ,
$$R\cong \Bbb Z_{p^2}[x,y_1,\cdots,y_m]/<\{g(x), px, y_sy_t,py_r,|1\le r,s,t\le
m\}
>$$
for some polynomials $g(x)\in \Bbb Z_{p^2}[x]$ such that
$\ol{g}(x)$ is irreducible over $\Bbb Z_p$, or $R\cong \Bbb
Z_{p^2}[x,y_1,\cdots,y_m]/M,$ for some
$$M=<\{g(x)^2,pg(x),g(x)y_r, px, y_sy_t,py_r,g(x)-\sum_{i=1}^mv_i(x)y_i|1\le r,s,t\le
m\}>$$ \nin where $g(x)$ is some polynomial over $\Bbb Z_{p^2}$
such that} $\mbox{deg}$\,$g(x)>1$ {\it and that $\ol{g}(x)$ is
irreducible over $\Bbb Z_p$, and at least one of the $v_i(x)$ is
not zero modulo $p$, while} $\mbox{deg}$\,$v_i(x)$ {\it is less
than }$\mbox{deg}$\,$g(x)$.

\vs{3mm}Finally, we remark that each ring of the four types in
Theorem 1 and Theorem 2 is obviously a $Z$-local ring. We also
remark that not all finite local rings whose zero-divisor graph is
uniquely determined are $Z$-local. For example, each of $\Bbb
Z_2[x_1,x_2,\cdots,x_n]/<x_1^2,x_2^2,\cdots,x_n^2>$ and $\Bbb
Z_4[x_1,x_2,\cdots,x_n]/<x_1^2,x_2^2,\cdots,x_n^2>$ is a finite
local rings with the property that $J(S)=Z(S)$ and $x^2=0, \forall
x\in Z(S)$. Obviously they are not $Z$-local.

\vs{5mm}\begin{center}{\bf REFERENCES}\end{center} \vs{3mm}

\nin 1.\, Anderson D.F.; Livingston P.S. The zero-divisor graph of
a commutative ring. J. of Algebra {\bf 1999},217, 434-447.

\nin 2.\, Lu Dancheng and Wu Tongsuo. The zero-divisor graphs
which are uniquely determined by neighborhoods. Preprint 2005.

\nin 2.\, Hideyuki Matrumura. Commutative Ring Theory. Cambridge
Studies in Advanced Mathematics, 1986.

\end{document}